\documentclass[12pt]{amsart}
\usepackage{amsfonts}
\usepackage{amsmath}
\usepackage{amssymb}
\usepackage{amsthm}
\usepackage{comment}
\usepackage{hyperref}
\usepackage{enumitem}

\makeatletter
\@namedef{subjclassname@2020}{%
  \textup{2020} Mathematics Subject Classification}
\makeatother

\newtheorem{thm}{Theorem}[section]
\newtheorem{lemma}[thm]{Lemma}
\newtheorem{prop}[thm]{Proposition}
\newtheorem{cor}[thm]{Corollary}
\newtheorem{example}[thm]{Example}

\newtheorem{mainthm}[thm]{Main Theorem}

\numberwithin{equation}{section}   

\theoremstyle{definition}
\newtheorem{defn}{Definition}[section]

\theoremstyle{remark}

\newtheorem{open}{Open question}

\newtheorem*{xrem}{Remark}

\begin{document}


\baselineskip=17pt


\title[Square Functions in $L^1$]{Square Functions for Ritt Operators in $L^1$}
 
\author[J. Hults]{Jennifer Hults}
\address{
Department of Mathematics and Statistics\\
University at Albany, SUNY, Albany, NY 12222 \\} 
\email{jhults@albany.edu} 

\author[K. Reinhold]{Karin Reinhold-Larsson } 
\address{
Department of Mathematics and Statistics\\
University at Albany, SUNY, Albany, NY 12222 } 
\email{reinhold@albany.edu} 

\date{}




\maketitle

\begin{abstract}
A power bounded operator $T$ is a Ritt operator in $L^p$ if $\sup_n n\lVert T^n-T^{n+1}\rVert<\infty$. When $T$ is a positive contraction and a Ritt operator in $L^p$, $1<p<\infty$, the square function \\
 $\Big( \sum_n n^{2m+1} |T^n(I-T)^{m+1}f|^2 \Big)^{1/2}$ is a bounded operator \cite{LeMX-Vq} in $L^p$. We show that if $T$ is a Ritt operator in $L^1$, then
\[Q_{\alpha,s,m}f=\Big( \sum_n n^{\alpha} |T^n(I-T)^mf|^s \Big)^{1/s}\] is bounded $L^1$ when $\alpha+1<sm$, and examine related questions regarding variational and oscillation norms.
\end{abstract}

\section{Introduction}

The purpose of this paper is to show that a certain "square" function is bounded in $L^1$ for a class of contractions operators.
Square functions and variational inequalities have been fruitful in martingale theory, harmonic analysis, in particular in Littlewood-Paley theory, and ergodic theory, to bound associated maximal functions, establish convergence of sequences of operators and control their speed of convergence. 


In the study of almost everywhere convergence of the powers $T^n f(x)$, Stein \cite{Stein} showed that for $T$ self adjoint with $\lVert T \rVert \le 1$, 
$Mf=\sup_{n\ge 1} |T^nf| $ and $Sf=\Big( \sum_{n\ge 1} n |T^{n+1}f-T^n f |^2 \Big)^{1/2}$ are bounded operators in $L^2$,
$\lVert Mf \rVert_2 \lesssim \lVert f \rVert_2$, and $\lVert Sf \rVert_2 \lesssim \lVert f \rVert_2.$

\begin{defn}\label{def:Stolz} A Stolz domain $\Sigma$ is the region within a unit disk characterized by $|1-z|\le c(1-|z|)$, for some constant $c$. 
$\Sigma_{\gamma}$ is the Stolz domain consiting of the interior of the convex hull of 1 and the disk $D(0,\sin \gamma)$.
\end{defn}
For $T$ normal $L^2$ contraction  with $\sigma(T)$ inside a Stolz domain  $\Sigma$,
Below, Jones and Rosenblatt  \cite{BJR-power}    
showed that both $Mf$ and $Sf$ are bounded operators in $L^2$; 
and if $T$ is an $L^1-L^{\infty}$ contraction, then $Mf$ is a bounded operator in $L^p$ for all $1<p<\infty$.

Around the time of the above results, the theory of functional calculus for sectorial operators was being develop, which provided control of square functions for a wider class of operators \cite{LeM-H,HH}.

\begin{defn} Let $Y$ be a Banach space. We say that  $T\in  \mathcal L(Y)$ is a Ritt operator if $T$ is power bounded and $\sup_n n\lVert T^n-T^{n+1}\rVert <\infty$.
\end{defn}

Ritt operators have been studied by many \cite{Ritt,Lyub,Nev,NZ,Bl}. 
The control of the resolvent of the operator is intricately linked to the containment of its spectrum within a Stolz domain.
Le Merdi \cite{LeM-H} exploited the connection between Ritt operators and sectorial operators to develop the $H^{\infty}$ functional calculus for Ritt operators, yielding new insights regarding the convergence of powers and maximal and square functions estimates \cite{LeM-H,LeMX-max,LeMX-Vq}. 

In particular,
LeMerdi and Xu \cite{LeMX-Vq} (Prop. 4.1, Thm. 4.4 \& Thm. 5.6) established the following implications for Ritt operators in $L^p$, $1<p<\infty$. 

\begin{thm} \label{thm:LeM} Let $(X,m)$ be a $\sigma$-finite measure space, $1<p<\infty$, and $T$ a positive contraction of $L^p(X,m)$. 
If $\sup_n n \lVert T^n-T^{n+1} \rVert<\infty$, then, for any fixed integer $m\ge 0$ and any real number $s>2$, 
\begin{itemize}
\item $\Big( \sum_n (n+1)^{2m+1} |T^n(I-T)^{m+1}f|^2 \Big)^{1/2}$, 
\item $\lVert \{T^nf\} \rVert_{v(s)}$, and more generally $\lVert \{n^m T^n(I-T)^mf\} \rVert_{v(s)}$, 
\item for any increasing sequence $\{n_k\}$, $\lVert \{T^nf\} \rVert_{o(2)}$, and more generally $\lVert \{n^m T^n(I-T)^mf\} \rVert_{o(2)}$, 
\end{itemize}
are bounded in $L^p$.
\end{thm}

Cohen, Cuny and Lin  \cite{CCL}   provided equivalent conditions between the Ritt, spectral conditions and the square functions.
\begin{thm}\label{thm:CCL}
Let $(X,m)$ be a $\sigma$-finite measure space, $1 <p <\infty$ and $T$ a positive contraction on $L^p(X, m)$. Then the following are equivalent:
\begin{enumerate}
\item  $\sup_n n\lVert T^n-T^{n+1}\rVert <\infty$,
\item there exists a constant $C_p>0$ such that $\lVert (\sum_n n |T^nf-T^{n+1}f|^2)^{1/2}\rVert \le C_p \lVert f\rVert_p$,
\item there exists a closed Stolz region $\Sigma$ and a constant $K>0$ such that
$\lVert u(T)\rVert \le K \sup_{z\in \Sigma} u(z)$ for every rational function $u$ with poles outside $\Sigma$.
\end{enumerate}
\end{thm}

Results in $L^1$ turned out to be more elusive even for the maximal function $Mf$. 
Let $(X,\beta,m)$ be a non-atomic, separable probability space and $\tau$ an invertible measure-preserving transformation on $(X,\beta,m)$. 
A probability measure $\mu$ on $\mathbb Z$ defines the contraction operator
$\tau_{\mu}f(x) = \sum_k \mu(k) f(\tau^k x)\label{eq:1}$, for $x\in X$, $f\in L^p(X)$, with $p\ge 1$. 
For such measures let $\hat{\mu}(t)=\sum_k \mu(k) e^{2\pi i kt}$ and $\tau_{\mu}$ is $\sigma(\mu)=\{\hat{\mu}(t): |t|\le 1/2\}$.  We say $\mu$ satisfies the bounded angular property iff $\sigma(\mu)\subset$ a Stoltz domain.

 Bellow and Calderon   \cite{BC}  showed that the maximal function \\
 $Mf=\sup_{n\ge 1} |\tau^n_{\mu}f|$ is weak (1,1) 
 for centered measures $\mu$ with finite second moment. Such measures have the bounded angular property.
 Losert \cite{Losert1,Losert2} constructed measures $\mu$ without the bounded angular property, for which pointwise convergence of $\tau^n_{\mu}$ failed.  Wedrychowicz  \cite{Chris} gave conditions for $Mf$ to be weak (1,1) for measures with $\sigma(\tau_{\mu})\subset$ a Stolz region but without a finite second moment.  

Regarding square functions in $L^1$, 
Jones, Ostrovskii and Rosenblat \cite{JOR} (Theorem 2.6) considered properties of square functions for differences of ergodic averages $A_nf(x) =\frac 1n \sum_{k=1}^n f(\tau^k x)$.
They showed that if $\{n_k\}$ is an increasing sequence of integers, then
$\Big( \sum_{n=1}^{\infty} \Big| A_{n_{k+1}}f-A_{n_k}f \Big|^2 \Big)^{1/2}$
is strong $(p,p)$ for $1<p\le 2$, and weak (1,1). It is also strong $(p,p)$ for $1<p<\infty$ if $\{n_k\}$ is lacunary.
And for sequences $n_k= p(k)$ for some
polynomial $p$ of degree $s > 0$, 
$\Big( \sum_{n=1}^{\infty} \sup_{n_k\le n \le n_{k+1}}  \Big| A_{n}f-A_{n_k}f \Big|^2 \Big)^{1/2}$
is strong $(p,p)$ for $1<p<\infty$ and weak (1,1). 


\begin{defn} 
Let $S=\{\nu_n\}$ a sequence of finite signed measures in $\mathbb Z$. 
Let $\Lambda$ be the semigroup generated by the support of the measures, $\Lambda=<\cup_n \text{supp}(\nu_n)>$. 
We say that $T\in \mathcal L(Y)$ is $S$-bounded (or $\Lambda$-power bounded) if $\lVert T\rVert_S=\sup_{n\in \Lambda} \lVert T^n\rVert <\infty$.
\end{defn}

With this notation, power bounded is equivalent with $\mathbb Z^+$-power bounded.

Note that, if at least one of the supports of the measures $\nu_n \in S$ is not contained in the positive integers, we are implicitly requiring the operator to be invertible.

When the sequence $S$ is the collection of convolution measures of a single measure $\mu^n=\mu*\mu*\cdots*\mu$, we denote it by $S_{\mu}=\{\mu^n\}_{n=1}^{\infty}$.

\begin{defn} 
 Let $S=\{\nu\}$ a single finite signed measure on the integers, and 
$T$ an S-bounded linear operator. 
Define the operator induced by $\nu$ as
\[T_{\nu}f=\sum_k \nu(k) T^kf.\] 
\end{defn}
With the restrictions on the operator $T$, $T_{\nu}$ is well defined for any $f\in Y$. 

For $0<\alpha<1$, consider the series expansion for 
\[(1-x)^{\alpha}=1-\sum_{k\ge 1} g(\alpha,k) x^k, \quad |x|\le 1.
\] 
The coefficients satisfy $g(\alpha,k)=\frac{\alpha |\alpha-1| \ldots |\alpha-k+1|}{k!} \ge 0$ 
and it is easy to check that $\sum_k g(\alpha,k)=1$ (\cite{Dunn}, \cite{LinDer}). Let
$\nu_{\alpha}$ 
be the probability measure on $\mathbb Z^+$ define by them,
 \begin{equation}\nu_{\alpha}(k)= g(\alpha,k) \label{eqn:valpha}.\end{equation}

\begin{defn} 
Let $T\in \mathcal L(Y)$ be power bounded.
Define $(I-T)^{\alpha}$ 
as \[(I-T)^{\alpha} = I - T_{\nu_{\alpha}}= I - \sum_{k\ge 1} g(\alpha,k) T^k.\]
\end{defn}

By separating the integer part from the fractional part, we can define  $(I-T)^m$ for any real $m>0$. 
In particular,
$T_{\mu}^n(I-T_{\mu})^m f = T_{\nu_{n,m}} f$
where $\nu_{n,m}$ is a signed measure on $\mathbb Z$ satisfying $\hat \nu_{n,m}(t) = \hat \mu^n(t) (1-\hat\mu(t))^m$.

\begin{defn} 
For real numbers $s,r>0$, $\alpha>-1$, define the "generalized square function" (associated with the operator $T$) as:
\[Q_{\alpha,s,r}f =Q^{T}_{\alpha,s,r}f = \Big(\sum_{n=1}^{\infty} n^{\alpha} |T^n(I-T)^r f|^{s} \Big)^{1/s}.\]
\end{defn}

\begin{thm} \label{thm:RittResult} Let $(X,m)$ be a $\sigma$-finite measure space and $T$ a Ritt operator in $L^1(X)$, $r>0$ and $s\ge 1$. Then 
\begin{enumerate}
\item $Q^T_{\alpha,s,r}f$ 
is bounded in $L^1$ for $sr>\alpha+1$;
\item $\sup_n n^{\alpha} |T^n (I- T)^{r} f|$ is bounded in $L^1$ for any  $\alpha<r$;
\item   both  $\lVert n^{\beta} T^n(I-T)^r f\rVert_{v(s)}$ and $\lVert n^{\beta} T^n(I-T)^r f\rVert_{o(s)}$ are bounded in $L^1$ if $\beta=0$ or $\beta>0$ and $s(r-\beta)>1$;
\item for any increasing sequence $\{n_k\}$ such that $ n_{k+1} -n_k \sim n_k^{\alpha}$ for some $\alpha\in(0,1)$, then 
\begin{enumerate}
\item 
$\Big( \sum_k n_k^{\beta s} |T^{n_k}f-T^{n_{k+1}}f|^s \Big)^{1/s}$ and \\
$\Big( \sum_k n_k^{\beta s}  \max_{ n_k\le n<n_{k+1} } |T^n f-T^{n_k} f|^s \Big)^{1/s} $  
are bounded in $L^1$ for $s>(1-\alpha)/(1-\alpha-\beta)$; and\\
\item 
$\Big( \sum_k n_k^{\beta s}  \lVert \{T^n f:n_k\le n<n_{k+1}\} \rVert_{v(s)}^s\Big)^{1/s}$
is bounded in $L^1$ for $s>1/(1-\alpha-\beta)$.
\end{enumerate}
\end{enumerate}
\end{thm}

This theorem is proven in stages by first considering results for the special case of $T_{\mu}$.

\begin{xrem}
Note that, by Theorem \ref{thm:CCL},  $Q_{1,2,1}f$  is bounded in $L^p$ for any $p>1$, but in $L^1$, 
$Q_{1,s,1}f$ is bounded for $s>2$. And in general, when $r=\alpha$, $Q_{r,s,r}f$ is bounded in $L^1$ for $s>1+1/r$.
Additionally, by Theorem \ref{thm:LeM}, $Q_{2r+1,2,(r+1)}f$ is bounded in $L^p$ for any $p>1$. In $L^1$,
 $Q_{2r+1,s,(r+1)}f$ bounded  for any $s>2$,  and $Q_{\alpha,2,(r+1)}f$ is bounded  for any $\alpha<2r+1$.
\end{xrem}

\begin{open} Is $Q_{\alpha,s,r}f$ weak (1,1)  for $sr=\alpha+1$? 
\end{open}

\section{The case of $T_{\mu}$}

 First we address sufficient conditions under which $T_{\mu}$ is a Ritt operator in $L^1$. 
 For $T$ an $L^1-L^{\infty}$ contraction, Blunke \cite{Bl} showed, with a clever interpolation argument, that if $T$ is Ritt in $L^2$, it is Ritt in $L^p$, $1<p<\infty$. The technique does not extend to $L^1$. But for the particular case of $T_{\mu}$,
 Dungey  \cite{Dunn} established $T_{\mu}$ is Ritt in $L^1$ for measures with a stronger property than bounded angular ratio.
 
 \begin{thm} (Dungey \cite{Dunn}, Theorem 4.1) \label{thm:D}
Let $\mu$ be a probability measure on $\mathbb Z^+$ for which there exists $0<\alpha,c<1$ such that (i) $\mbox{Re} \, \hat\mu(t) \le 1-c |t|^{\alpha}$, and 
(ii)  $|\hat\mu'(t)|\lesssim |t|^{\alpha-1}$ for $0<|t|\le 1/2$.
If $T \in  \mathcal L(Y) $ is power bounded, then $T_{\mu}$ is Ritt in $L^1$.
\end{thm}

 Inspired by \cite{Chris} and 
 \cite{Dunn}, Cuny  \cite{Cuny-weak} considered more general spectral conditions for measures with support on $\mathbb Z$.

\begin{defn} 
Let $\mu$ be a probability measure on $\mathbb Z$ with Fourier transform $\hat\mu(t)=\sum_k \mu(k) e^{-2\pi i kt}$ that
is twice continuously differentiable on $0<|t|<1$. When there exists a continuous function $h(t)$ on $|t|\le 1$ with $h(0)=0$, $h(-t)=h(t)$, continuously differentiable on $0<|t|<1$ satisfying the following conditions:
(i)  $|\hat\mu(t)|\le 1-c\, h(t)$,
(ii)  $|t \, \hat\mu'(t)|\le c\, h(t)$,
(iii)  $|\hat\mu'(t)|\le c\, h'(t)$, and
(iv)  $|t \, \hat\mu''(t)|\le c\, h'(t)$,
 we say $\mu$ satisfies condition M1. 
 
 If in addition $\hat\mu$ satisfies
(v) $h(t)\le c th'(t)$ for $0<t<1$, then we say $\mu$ satisfies condition M, or $\mu$ is an M-measure.
\end{defn}

Note that  M1 implies bounded angular ratio
\[ |1-\hat\mu(t)| \le | \int_0^t \hat\mu'(u) du | \le \int_0^{|t|} h'(u) du = h(t) \le \frac{1-|\hat\mu(t)|}{c}.\]

If $\mu$ is a centered measure with finite second moment, then it satisfies M with $h(t)=t^2$. 
The next example exhibits a non-centered measure  without finite first moment.
\begin{example}   For fixed $0<\alpha<1$, let $\nu_{\alpha}$ be the probability measure on $\mathbb Z$ defined in (\ref{eqn:valpha}). 
Then $\hat\nu_{\alpha}(t) = 1-(1-e^{2\pi i t})^{\alpha}$. From \cite{Dunn}, $\nu_{\alpha}$ satisfies the conditions of Theorem \ref{thm:D} and hence $T_{\nu_{\alpha}} = I-(I-T)^{\alpha}$ is a Ritt operator in $L^1$. Additionally, $\nu_{\alpha}$ satisfies property M
 (\cite{Cuny-weak} (Proposition 3.3)). 
 
\end{example}

More examples of measures satisfying M can be found in \cite{Cuny-weak}.

\begin{thm} \label{thm:weak} 
\cite{Cuny-weak} 
Let $(X,m)$ be a probability space, $\tau$ a measure preserving transformation on it.
Let $\mu$ be a probability measure on $\mathbb Z$ with property M1. Then 
$m(x\in X: \sup_n |\tau_{\mu}^nf(x)|>\lambda) \le \frac{C}{\lambda} \lVert f\rVert_1$
and $\tau_{\mu}$ is weak-$L^1$-Ritt:  \\
$ m(x\in X: \sup_n n^r |(\tau_{\mu}^n-\tau_{\mu}^{n+r})f(x)|>\lambda) \le \frac{C_r}{\lambda} \lVert f\rVert_1$
for any $f\in L^1(X)$.\\
If in addition, $\mu$ is an M-measure, then $\tau_{\mu}$ is Ritt in $L^1$.
\end{thm} 

Cuny proved this result for the  shift operator in $l^1(\mathbb Z)$ but a transfer argument extends his result to operators defined by a measure presrving transformation on a probability space.


\begin{mainthm}\label{thm:main}
Let $\mu$ be an M-probability measure on $\mathbb Z$ and $T\in \mathcal L(X)$ be $S_{\mu}$-bounded.
 If $r>0$, and  $sr>\alpha+1$,
then $Q_{\alpha,s,r}^{T_{\mu}}f$ is a bounded operator in $L^1$.
\end{mainthm}

\begin{cor} \label{cor:sup} Let $\mu$ be an M-probability measure on $\mathbb Z$ and $T\in \mathcal L(X)$ be $S_{\mu}$-bounded.
For $0\le \alpha<r$,
$\sup_n n^{\alpha} |T_{\mu}^n (I- T_{\mu})^{r}f|$ is bounded in $L^1$ and 
$\lim_{n\to \infty} n^{\alpha} |T_{\mu}^n (I- T_{\mu})^{r}f| = 0$ in norm and a.e..
\end{cor}

Theorem \ref{thm:LeM} showed that $\lVert n^{r} T_{\mu}^n(I-T_{\mu})^rf\rVert_{v(s)}$ and $\lVert n^{r} T_{\mu}^n(I-T_{\mu})^rf \rVert_{o(s)}$ 
are bounded in $L^p$, for $s>2$ and $p>1$ , even in the case $r=0$.
In $L^1$, we obtained the following variations and oscillations results.
\begin{prop} \label{prop:genvar}
Let $\mu$ be an M-probability measure on $\mathbb Z$ and $T\in \mathcal L(X)$ be $S_{\mu}$-bounded.
Let $s\ge 1$, $r>0,$ and  $\beta\ge 0$  be fixed. If  $\beta=0$ or if $s(r-\beta)>1$, then both 
 $\lVert n^{\beta} T_{\mu}^n(I-T_{\mu})^r f \rVert_{v(s)}$ 
 and
 $\lVert n^{\beta} T_{\mu}^n(I-T_{\mu})^r f \rVert_{o(s)}$ are bounded in $L^1$.
\end{prop}

\begin{open} Are there values of $\beta>-1$ and $s>0$ for which
$\lVert n^{\beta} T_{\mu}^nf\rVert_{v(s)}$ and $\lVert n^{\beta} T_{\mu}^nf \rVert_{o(s)}$ are bounded in $L^1$.
\end{open}

The next result handles cases of 
differences along subsequences with increasing gaps.

\begin{prop} \label{prop:longvar}
Let $\mu$ be an M-probability measure on $\mathbb Z$ and $T\in \mathcal L(X)$ be $S_{\mu}$-bounded.
Let  $\{n_k\}$ be an increasing sequence such that $ n_{k+1} -n_k \sim n_k^{\alpha}$ 
for some $0<\alpha<1$, and $0\le \beta<1-\alpha$. Then
\begin{enumerate}
\item 
$\Big( \sum_k n_k^{\beta s} |T_{\mu}^{n_k}f-T_{\mu}^{n_{k+1}}f|^s \Big)^{1/s}$
is bounded in $L^1$ for $s>(1-\alpha)/(1-\alpha-\beta)$.
 In particular,  $\Big( \sum_k  |T_{\mu}^{n_k}f-T_{\mu}^{n_{k+1}}f|^s \Big)^{1/s}$ is bounded in $L^1$ for $s>1$;
\item 
$\Big( \sum_k n_k^{\beta s}  \lVert \{T_{\mu}^n f:n_k\le n<n_{k+1}\} \rVert_{v(s)}^s\Big)^{1/s}$
is bounded in $L^1$ for $s>1/(1-\alpha-\beta)$; \\
and $\Big( \sum_k n_k^{\beta s}  \max_{ n_k\le n<n_{k+1} } |T_{\mu}^n f-T_{\mu}^{n_k} f|^s \Big)^{1/s} $  is
bounded in $L^1$ for $s>(1-\alpha)/(1-\alpha-\beta)$.
\end{enumerate}
\end{prop}


\section{Auxiliary Lemma}
In these notes, $c$ and $C$ denote constants whose values may change from one instance to the next. We simplify $e(x)=e^{2\pi i x}$. And for  $0\le x,y$, we say $x \lesssim y$ if there exists a constant $c>0$ such that $x \le c y$.

\begin{lemma}\label{lem:basic} 
Let $S=\{\Delta_n\}$ be a sequences of (finite) signed measures on $\mathbb Z$ and $T\in \mathcal L(Y)$ an $S$-bounded operator.
Let
\begin{align*} 
A= & \int_{|t|<1/2} \frac 1{|t|} \Big( \sum_n |\hat \Delta_n(t)|^s \Big)^{1/s} dt,  \quad
B= & \int_{|t|<1/2} |t| \Big( \sum_n |\hat \Delta''_n(t)|^s \Big)^{1/s} dt, \\
 C=& \sum_{k\neq 0} \frac 1{|k|} \Big( \sum_n |\hat \Delta_n(1/|k|)|^s \Big)^{1/s},  \quad
 D=& \sum_{k\neq 0} \frac 1{|k|^2} \Big( \sum_n |\hat \Delta'_n(1/|k|)|^s \Big)^{1/s}, \\
 E= &\Big( \sum_n  \Big|   \Delta_n(0)   \Big|^{s} \Big)^{1/s}. \quad &
\end{align*}
If $A, B, C, D$ and $E$ are all finite,  then, for any $f\in X$,
\[ \Big\lVert  \Big( \sum_n |T_{\Delta_n}f|^s \Big)^{1/s} \Big\rVert \lesssim \lVert f\rVert.\]
\end{lemma}

\begin{proof}
Without loss of generality, we assume $\lVert T \rVert_S=1$.

\begin{align*}   \Big\lVert \Big( \sum_n |T_{\Delta_n}f|^s \Big)^{1/s} \Big\rVert
= & \Big\lVert \Big(\sum_n  \Big| \sum_k \Delta_n(k) T^kf\Big|^{s} \Big)^{1/s} \Big\rVert \\
\le &  c \Big\lVert   \Big( \sum_n  \Big|  \sum_{k\neq 0} \int_{|t|<1/2|k| } \hat \Delta_n(t) e(kt) dt  \,  T^kf \Big|^{s} \Big)^{1/s} \Big\rVert \\
 & + c \Big\lVert    \Big( \sum_n  \Big|  \sum_{k\neq 0}   \int_{1/2|k|<|t|<1/2} \hat \Delta_n(t) e(kt) dt  \,  T^kf \Big|^{s} \Big)^{1/s} \Big\rVert \\
& + c \Big\lVert    \Big( \sum_n  \Big|   \Delta_n(0)   \Big|^{s} \Big)^{1/s} |f| \Big\rVert \\
 =  & c (  \mbox{I} + \mbox{II}+  E \lVert f \rVert).
\end{align*}

For the first term we have,
\begin{align*}
  \mbox{I} 
=  & \Big\lVert   \Big( \sum_n  \Big|  \sum_{k\neq 0} \int_{|t|<1/2|k|} \hat \Delta_n(t) e(kt) dt  \,  T^kf\Big|^{s} \Big)^{1/s} \Big\rVert \\
\le & \lVert f\rVert \sum_{k\neq 0} \Big( \sum_n \Big|   \int_{|t|<1/ 2|k| } \hat \Delta_n(t) e(kt) dt   \Big|^{s} \Big)^{1/s}\\
\le &   \lVert f\rVert \int_{|t|<1/2} \frac1{|t|}  \Big( \sum_n  \Big| \hat \Delta_n(t) \Big|^s \Big)^{1/s}  \, dt
= A \lVert f\rVert.
  \end{align*}

For the second term we have,
\begin{align*}
\mbox{II} \le & \Big\lVert   \Big( \sum_n  \Big|  \sum_{k\neq 0}   \int_{1/2|k| <|t|<1/2} \hat \Delta_n(t) e(kt) dt  \,  T^kf \Big|^{s} \Big)^{1/s} \Big\rVert \\
\le & \Big\lVert  \sum_{k\neq 0}  \Big( \sum_n  \Big |  \int_{1/2|k|<|t|<1/2} \hat \Delta_n(t) e(kt) dt \Big |^s \Big)^{1/s} |T^kf | \Big\rVert\\
\le & \lVert f\rVert  \sum_{k\neq 0}  \Big( \sum_n  \Big |  \int_{1/2|k|<|t|<1/2} \hat \Delta_n(t) e(kt) dt \Big |^s \Big)^{1/s}.
\end{align*}

The integrand decomposes as
\begin{align*}
\Big| \int_{1/|k| <|t|<1/2} \hat \Delta_n(t) e(kt) dt \Big| \le & 
 \Big|\int_{1/|k| <|t|<1/2} \hat \Delta'_n(t) \frac{e(kt)}{2\pi k} dt \Big| \\
& + \Big|\frac{\hat \Delta_n(1/|k| ) e(k/|k| )}{2\pi k} - \frac{\hat \Delta_n(-1/|k| ) e(-k/|k| )}{2\pi k} \Big|\\  
\le &  \Big| \int_{1/|k| <|t|<1/2} \hat \Delta''_n(t) \frac{e(kt)}{4\pi^2 k^2} dt \Big| \\
 & + \Big| \frac{\hat \Delta_n(1/|k| ) e(k/|k| )}{2\pi k} - \frac{\hat \Delta_n(-1/|k| ) e(-k/|k| )}{2\pi k} \Big|\\  
 &  +  \Big| \frac{\hat \Delta'_n(1/|k| ) e(k/|k| )}{4\pi^2 k^2} - \frac{\hat \Delta'_n(-1/|k| ) e(-k/|k| )}{4\pi^2 k^2} \Big| 
 \end{align*}

\begin{align*}
\sum_{k\neq 0}  & \Big( \sum_n  \Big |  \int_{1/2|k|<|t|<1/2} \hat \Delta''_n(t) \frac{e(kt)}{4\pi^2 k^2}  dt \Big |^s \Big)^{1/s} \\
\lesssim  & \sum_{k\neq 0} \frac 1{k^2}  \int_{1/2|k|<|t|<1/2} \Big( \sum_n  \Big | \hat \Delta''_n(t) \Big |^s \Big)^{1/s} \ dt  \\
\lesssim & \int_{0<|t|<1/2}  |t| \,  \Big( \sum_n  \Big | \hat \Delta''_n(t) \Big |^s \Big)^{1/s}\, dt  = B .
\end{align*}

Thus, $\mbox{II}\lesssim (B+C+D) \lVert f\rVert$.

\end{proof}

The proof of this lemma can be adapted for the following setting.

\begin{lemma}\label{lem:basic2} 
Let $\Delta_n$ and $T$ as in Lemma \ref{lem:basic}, $\{n_k\}$ an increasing sequence and $I_k=[n_k,n_{k+1})$.
Let
\begin{align*}
A= & \int_{|t|<1/2} \frac 1{|t|}        \Big( \sum_k n_k^{\beta} \max_{n\in I_k} |\hat \Delta_{n}(t)|^s \Big)^{1/s} dt,  \quad
B= & \int_{|t|<1/2} |t|     \Big( \sum_k n_k^{\beta} \max_{n\in I_k} |\hat \Delta''_{n}(t)|^s \Big)^{1/s} dt, \\
C=&\sum_{l\neq 0} \frac 1{|l|}       \Big( \sum_k n_k^{\beta} \max_{n\in I_k} |\hat \Delta_{n}(1/|l|)|^s \Big)^{1/s},  \quad
D= &\sum_{l\neq 0} \frac 1{|l|^2}       \Big( \sum_k n_k^{\beta} \max_{n\in I_k} |\hat \Delta'_{n}(1/|l|)|^s \Big)^{1/s}, \\
 E= &     \Big( \sum_k n_k^{\beta} \max_{n\in I_k}  \Big|   \Delta_{n}(0)   \Big|^{s} \Big)^{1/s}. &
\end{align*}
If $A, B, C, D, E$ are all finite, then, for any $f\in X$,
\[ \Big\lVert    \Big( \sum_k n_k^{\beta} \max_{n\in I_k} |T_{\Delta_{n}}f|^s \Big)^{1/s} \Big\rVert \lesssim  \lVert f \rVert.\]
\end{lemma}

\section{Proofs of Results}

\noindent {\textit Proof of Theorem \ref{thm:main}:}\\
Let $\Delta_n$ be
the measure on the integers defined by $\hat\Delta_n=n^{\alpha/s} \hat\mu^n (1-\hat\mu)^r$, that is, 
$T_{\Delta_n}=n^{\alpha/s} T_{\mu}^n (I-T_{\mu})^r$. Using Lemma \ref{lem:basic} with $\Delta_n$, it suffices to show that the corresponding terms defined in the lemma, A, B, C, D and E, are bounded.
The case for $\alpha<-1$ is immediate. Assume $\alpha\ge -1$. 

For $\alpha>-1$, dominating the sum by an integral and using condition (i) of property M, we obtain
\[\sum_n n^{\alpha}  |\hat\mu(t)|^{ns} \le \sum_n n^{\alpha}  (1-c \, h(t))^{ns}    \lesssim   \frac 1{h(t)^{\alpha+1}}.\]
For $\alpha=-1$, 
\[\sum_n \frac 1n  |\hat\mu(t)|^{ns} \le \sum_n \frac 1n  (1-ch(t))^{ns}= | \ln(1-(1-ch(t))^s) | 
 \lesssim \frac 1{h(t)^{\gamma}}
\mbox{ for any }\gamma>0.\]

Condition (iii) of property M also implies
$|1-\hat\mu(t)|  \le c \int_0^{|t|} h(s)' ds = c h(t)$. 

For term A, we estimate
\begin{align*}
A=\int_{|t|<1/2} \frac1{|t|} &  \Big( \sum_n  \Big| \hat\Delta_n(t) \Big|^s \Big)^{1/s}  \\
= &    \int_{|t|<1/2}  \Big( \sum_n n^{\alpha} | \hat\mu(t)|^{ns} \Big)^{1/s} \frac{|1-\hat\mu(t)|^{r}}{|t|}  \\
\lesssim  & \begin{cases} \int_{0<t<1/2} h(t)^{r-(\alpha+1)/s-1}  \, h'(t) \, dt & \mbox{ for } \alpha>-1 \\
 \int_{0<t<1/2} h(t)^{r-\gamma/s-1}  \, h'(t) \, dt & \mbox{ for } \alpha=-1.
 \end{cases}
  \end{align*}
When $\alpha>-1$, the integral is finite for $sr >\alpha+1$, and when $\alpha=-1$, the integral is finite for $r>0$ 
since we can choose $\gamma$ arbitrarily small, say $\gamma=sr/2$.

Similarly for E, with $sr >\alpha+1$ for $\alpha>-1$ and $0<\gamma<sr /2$ for $\alpha=-1$,
\begin{align*}
E= \Big( \sum_n  \Big|   \Delta_n(0)   \Big|^{s} \Big)^{1/s} \le &
\Big( \sum_n n^{\alpha}    \int_{|t|<1/2} |\hat\mu(t)|^{ns} |1-\hat\mu(t)|^{sr}\, dt \Big)^{1/s} \\
\lesssim  &\begin{cases} \Big(  \int_{|t|<1/2}  h(t)^{ (sr-(\alpha+1))}\, dt \Big)^{1/s} <\infty, & \mbox{ for } \alpha>-1 \\
 \Big(  \int_{|t|<1/2}  h(t)^{ (sr-\gamma)} \, dt \Big)^{1/s} <\infty, & \mbox{ for } \alpha=-1. 
  \end{cases}
\end{align*}

With $B_n=T^n (I-T)^r$,
\[
|\hat B'_n(t)| = | n \hat\mu^{n-1}(t) (\hat \mu'(t)) (1-\hat\mu(t))^r 
     - \hat\mu^{n}(t) \hat \mu'(t)(1-\hat\mu(t))^{r-1}|, \]
and 
\begin{align*}
|\hat B''_n(t)| = &| n(n-1) \hat\mu^{n-2}(t) (\hat \mu'(t))^2 (1-\hat\mu(t))^r 
     +n \hat\mu^{n-1}(t) \hat \mu''(t)(1-\hat\mu(t))^r \\
& - 2n r \hat\mu^{n-1}(t) (\hat \mu'(t))^2 (1-\hat\mu(t))^{r-1}  \\
&+ r (r-1) \hat\mu^{n}(t) (1-\hat\mu(t))^{r-2}  (\hat\mu'(t))^2\\
& - r \hat\mu^{n}(t) (1-\hat\mu(t))^{r-1}  \hat \mu''(t)| \\
\lesssim & n^2 |\hat\mu^{n-2}(t)| \frac{h'(t)}{|t|} h^{r+1}(t) + 
n |\hat\mu^{n-1}(t)| \frac{h'(t)}{|t|} h^{r}(t) \\
& + |\hat\mu^{n}(t)| \frac{h'(t)}{|t|} h^{r-1}(t).
\end{align*}
Thus, for $\alpha>-1$,
\[
\Big(\sum_n   |\hat \Delta''_n (t)|^s=  \sum_n n^{\alpha}  |\hat B''_n (t)|^s\Big)^{1/s}  
\lesssim    
h(t)^{(r-1)-(\alpha+1)/s} \frac{|h'(t)|}{|t|}.\]
Since $(1+\alpha)< s r $,
\begin{align*}
B=   \int_{|t|<1/2} |t| \, \Big( \sum_n   |\hat\Delta''_n(t)|^s  \Big)^{1/s} dt   
 \lesssim &  \int_{0<t<1/2}   h(t)^{(r-1)-(\alpha+1)/s } \, h'(t) \, dt <\infty.   \label{eq:II1} 
\end{align*}

When $\alpha=-1$, choosing $0<\gamma<r$, the estimate is
\[   \int_{|t|<1/2} |t| \, \Big( \sum_n   |\hat\Delta''_n(t)|^s  \Big)^{1/s} dt   
 \lesssim   \int_{0<t<1/2}   h(t)^{r-\gamma-1 } \, h'(t) \, dt <\infty.  \]

For the remaining terms, we address the case $\alpha>-1$ since the estimates for the case $\alpha=-1$ follow similar arguments.
For $(1+\alpha)< s r $, we have
\begin{align*}
 C= \sum_{k\neq 0} \frac 1{|k|} &
\Big(    \sum_n   |\hat\Delta_n(1/|k| )|^s  \Big)^{1/s}  \\
\lesssim   & \sum_{k\neq 0} \frac 1{|k|} \Big(    \sum_n n^{\alpha} |\hat\mu(1/|k|)|^{ns} 
|1-\hat\mu(1/|k|)|^{s r}  \Big)^{1/s} \nonumber  \\
\lesssim & \sum_{k> 0} \frac 1k  h(1/k)^{r - (\alpha+1)/s}
\le c+ \int_0^{1/2} \frac{h(t)^{r - (\alpha+1)/s}}t \, dt <\infty.
\label{eq:II2}
\end{align*}
Since
\[ \Big( \sum_n  |\hat\Delta'_n(1/|k| )|^s\Big)^{1/s} = \Big( \sum_n n^{\alpha}  |\hat B'_n(1/|k| )|^s\Big)^{1/s}
 \lesssim  h(1/k)^{r-(\alpha+1)/s} h'(1/k), \]
 it follows 
\begin{align*}
D = & \sum_{k\neq 0} \frac 1{|k|^2}
\Big(    \sum_n   |\hat\Delta'_n(1/|k| )|^s  \Big)^{1/s} \nonumber \\
 & 
  \lesssim     \sum_{k\neq 0} \frac 1{|k|^2}   h(1/k)^{r - (\alpha+1)/s} h'(1/k) <\infty.   \label{eq:II3}
\end{align*}~\hfill$\square$
  \bigskip
  
  \noindent {\textit Proof of Corollay \ref{cor:sup}:}\\
An application of Abel's summation yields
\begin{align*}
n^{\alpha } |T_{\mu}^{n-1}(I-T_{\mu})^{r} f| \le & \sum_{k=0}^{n-1} ((k+1)^{\alpha }-k^{\alpha } )|T_{\mu}^k(I - T_{\mu})^{r} f| \\
& + \sum_{k=1}^n k^{\alpha } |T_{\mu}^{k-1}(1-T_{\mu})^{r+1}f|\\
\lesssim & \sum_{k=0}^{n-1} (k+1)^{\alpha  -1} |T_{\mu}^k (I- T_{\mu})^{r} f| \\
& + \sum_{k=1}^n k^{\alpha } |T_{\mu}^{k-1}(1-T_{\mu})^{r+1}f|\\
\lesssim &  \ Q_{\alpha -1,1,r}f + Q_{\alpha ,1,r+1}f.
\end{align*}
By Theorem \ref{thm:main},
 both generalized square functions on the right are bounded in $L^1$ for $\alpha< r$ and $r>0$. Therefore 
$\sup_n n^{\alpha } |T_{\mu}^n(I - T_{\mu})^{r})f|$ is also bounded in $L^1$.
~\hfill$\square$ \bigskip

\noindent {\textit Proof of Proposition \ref{prop:genvar}:}\\

Let $\Delta_{n,r}f=T_{\mu}^n(I-T_{\mu})^r f$, and 
let $\{n_k\}$ any increasing sequence. 
\[D_{k,\beta}=  n_k^{\beta} \Delta_{n_k,r} - n_{k+1}^{\beta} \Delta_{n_{k+1},r} 
=  n_k^{\beta} \sum_{u=n_k}^{n_{k+1}-1} \Delta_{u,r+1} - (n_{k+1}^{\beta}-n_k^{\beta} ) \Delta_{n_{k+1},r}.
\]
\begin{align*}
\Big( \sum_k |D_{k,\beta}f|^s \Big)^{1/s}
\lesssim &  \sum_k  \sum_{u=n_k}^{n_{k+1}-1}  u^{\beta } |\Delta_{u,r+1}f|  
+ \Big( \sum_k n_k^{s\beta}  |\Delta_{n_{k},r}f|^s  \Big)^{1/s}.
\end{align*}
Thus, 
\[ \lVert n^{\beta}\Delta_{n,r}f \rVert_{v(s)}  \lesssim Q_{\beta ,1,r+1}f + Q_{\beta s,s,r}f,\]
which, by Theorem \ref{thm:main}, are bounded in $L^1$ for $s(r-\beta)>1$.

For $n_k\le n \le n_{k+1}$
\[ \max_{n_k\le n \le n_{k+1}}
\Big| n^{\beta} \Delta_{n,r} - n_{k}^{\beta} \Delta_{n_{k},r}f \Big|^s \le 
 \sum_{u=n_k}^{n_{k+1}-1} u^{\beta } |\Delta_{u,r+1}f| 
+  \sum_{u=n_k+1}^{n_{k+1}} n^{s\beta}  |\Delta_{n,r}f|^s.
\]
Thus, 
\[ \lVert n^{\beta}\Delta_{n,r}f \rVert_{o(s)}  \lesssim Q_{\beta ,1,r+1}f + Q_{\beta s,s,r}f\]
is bounded in $L^1$ for $s(r-\beta)>1$.

When $\beta=0$,
\[ \lVert \Delta_{n,r}f \rVert_{v(s)}  \lesssim Q_{0,1,r+1}f \quad
\mbox{ and }\quad
 \lVert \Delta_{n,r}f \rVert_{o(s)}  \lesssim Q_{0,1,r+1}f,\]
are bounded in $L^1$ for $r> 0$. 
\hfill$\square$

\bigskip
\noindent {\textit Proof of Proposition \ref{prop:longvar}:}  
Assume $\mu$ satisfy condition M and 
apply Lemma \ref{lem:basic} with $\Delta_k = n_k^{\beta } (T_{\mu}^{n_k} f - T_{\mu}^{n_{k+1}})$.

For $\gamma>\alpha$,
 \[\sum_k n_k^{s\gamma} (n_{k+1}-n_k) |\hat\mu(t)|^{n_k s}\lesssim \frac1{h(t)^{s\gamma+1}} .\]

\begin{align*}
 \mbox{A} = &
\int_{|t|<1/2} \frac 1{|t|} \Big( \sum_k  n_k^{\beta s}  |\hat\mu(t)|^{sn_k} \, |1-\hat\mu(t)^{(n_{k+1}-n_k)}|^s \Big)^{1/s} dt\\
\lesssim & \int_{|t|<1/2} \frac {|1-\hat\mu(t)|}{|t|} \Big( \sum_k n_k^{\beta s} (n_{k+1}-n_k)^s |\hat\mu(t)|^{s n_k}  \Big)^{1/s} \,  dt\\
\lesssim & \int_{|t|<1/2} \frac {|1-\hat\mu(t)|}{|t|} \Big( \sum_k n_k^{(\alpha+\beta-\alpha/s) s} (n_{k+1}-n_k) |\hat\mu(t)|^{s n_k}  \Big)^{1/s} \,  dt\\
\lesssim & \int_{|t|<1/2} \frac {h(t)}{|t|} \Big(\frac 1{h(t)^{s(\alpha+\beta-\alpha/s)+1}}\Big)^{1/s} \, dt \\
 \lesssim & \int_{|t|<1/2} \frac {h'(t)}{h(t)^{(\alpha+\beta)+(1-\alpha)/s}} \, dt   <\infty,
  \end{align*}
for $\alpha+\beta+(1-\alpha)/s<1, s>1$. 
Similarly, 
\begin{align*}
\mbox{E} \le &  \Big(  \sum_k n_k^{(\alpha+\beta)s}   \int_{|t|<1/2} |\hat\mu(t)|^{n_k s}  |1-\hat\mu(t)|^s dt    \Big)^{1/s} \\
\lesssim &  \Big(  \int_{|t|<1/2}   h(t)^{s(1-(\alpha+\beta)-(1-\alpha)/s)} dt   \Big)^{1/s}<\infty.
\end{align*}

For the next term,
\begin{align*}
n_k^{-\beta}  \hat\Delta''_k(t) = &  n_k(n_k-1) \hat\mu^{n_k-2}(t) (\hat \mu'(t))^2 - n_{k+1}(n_{k+1}-1) 
         \hat\mu^{n_{k+1}-2}(t)(\hat \mu'(t))^2 \\
&  +n_k \hat\mu^{n_k-1}(t) \hat \mu''(t) - n_{k+1} \hat\mu^{n_{k+1}-1}(t) \hat \mu''(t)  \\
= & n_k(n_k-1)\hat\mu^{n_k-2}(t) (1-\hat\mu^{n_{k+1}-n_k}(t)) (\hat \mu'(t))^2 \\
& +[n_k(n_k-1)-n_{k+1} (n_{k+1}-1)] \hat\mu^{n_{k+1}} (\hat \mu'(t))^2 \\
& +  n_k  \hat\mu^{n_k-1}(t) (1-\hat\mu^{n_{k+1}-n_k}(t)) \hat \mu''(t) \\
& +(n_k-n_{k+1}) \hat\mu^{n_{k+1}-1}(t) \hat \mu''(t).
\end{align*}
Estimating
\begin{align*}
n_k(n_k-1)-n_{k+1}(n_{k+1}-1)= & n_k^2-n_{k+1}^2+(n_{k+1}-n_k) \\
\sim & n_k^{\alpha}(n_{k+1}+n_k) + n_k^{\alpha} \\
\le  & n_k^{\alpha}(2n_k+n_k^{\alpha})  \lesssim n_k^{1+\alpha},
\end{align*}
we have
\begin{align*}
|\hat\Delta''_k(t)| \lesssim & \Big[   n_k^{2+\alpha+\beta}  \hat\mu^{n_k-2}(t) h(t)^2 
     + n_k^{1+\alpha+\beta} |\hat\mu(t)|^{n_{k+1}-1}  h(t)  \Big. \\
 & \Big. + n_{k}^{1+\alpha+\beta} |\hat\mu(t)|^{n_{k+1}-2}  h(t) 
 +n_{k}^{\alpha+\beta} |\hat\mu(t)|^{n_{k+1}-1}  \Big] \frac{h'(t)}{|t|}.
\end{align*}
Thus
\begin{align*}
\Big(   \sum_k   |\hat\Delta''_k(t)|^s    \Big)^{1/s} 
  \lesssim    \frac{1}{h(t)^{(\alpha+\beta)+(1-\alpha)/s) } }   \frac{h'(t)} {|t|}.
\end{align*}
%
%
Then,  for $\alpha+\beta+(1-\alpha)/s<1$,
\begin{equation*}
\mbox{B} =  \int_{|t|<1/2} |t|\,  \Big( \sum_k  |\hat\Delta''_{k}(t)|^s  \Big)^{1/s} \, dt 
  \lesssim   \int_{0<t<1/2} 
\frac {h'(t)}{h(t)^{(\alpha +\beta)+(1-\alpha)/s)) }}\, dt<\infty,
 \label{eq:II4}
\end{equation*}
and
\begin{align*}
\mbox{C} =&  \sum_{l\neq 0} \frac 1{|l|}
\Big(    \sum_k |\hat\Delta_{k}(1/|l| )|^s  \Big)^{1/s} \\
\lesssim  & \sum_{l\neq 0} \frac 1{|l|} \Big(    \sum_k  n_k^{(\alpha+\beta) s} |\hat\mu(1/|l|)|^{n_k s} 
|1-\hat\mu(1/|l|)|^s  \Big)^{1/s} \\
\lesssim  &  \sum_{l\neq 0} \frac {h(1/|l|)^{1-(\alpha+\beta)-(1-\alpha)/s)}}{|l|}   \label{eq:II5} 
<\infty.
\end{align*}

For the last term, we have
\begin{align*}
n_k^{-\beta} |\hat\Delta'_{k}(1/|l| )| \le &
 n_k | \hat \mu^{n_k-1}(1/|l|)| \,  
\Big| 1-\hat\mu(1/|l|)^{n_{k+1}-n_k} \Big| \,|\hat\mu'(1/|l|) | \\
& +  (n_{k+1}-n_k)  \Big| \hat\mu(1/|l|) \Big|^{n_{k+1}-1}  \, |\hat\mu'(1/|l|)|\\
 \lesssim & n_k^{1+\alpha} | \hat \mu^{n_k-1}(1/|l|)| \, |\hat\mu'(1/|l|)| |1-\hat\mu(1/|l|)| \\
& + n_k^{\alpha} | \hat \mu^{n_{k+1}-1}(1/|l|)| \, |\hat\mu'(1/|l|)|.
 \end{align*}
Then
\[ \Big( \sum_k  |\hat\Delta'_{k}(1/|l| )|^s\Big)^{1/s}
 \lesssim  \frac {h'(1/|l|)}{h(1/l)^{(\alpha+\beta)+(1-\alpha)/s )}} , \]
 and
\[
\mbox{D} \lesssim        \sum_{l>0} \frac 1{l^2} \frac {h'(1/|l|)}{h(1/l)^{(\alpha+\beta)+(1-\alpha)/s }} <\infty
\]
for $0<\alpha+\beta+(1-\alpha)/s<1$. 

 \bigskip

Now, let $I_k=[n_k,n_{k+1})$. 
\begin{align*}
\sum_k  n_k^{\beta s}  \lVert\{T_{\mu}^n f-T_{\mu}^{n_k}f: &n\in I_k \} \rVert_{v(s)}^s     \le   c \sum_k  n_k^{\beta s} \sum_{n\in I_k} |T_{\mu}^n f-T_{\mu}^{n_k}f|^s \\
\le  &c \sum_k n_k^{\beta s}  (n_{k+1}-n_k) \max_{n\in I_k} |T_{\mu}^n f-T_{\mu}^{n_k}f|^s .  
 \end{align*}
 Using Lemma \ref{lem:basic2} and arguments similar to the above, 
 \[\Big(\sum_k  n_k^{\beta s}  \lVert \{T_{\mu}^n f-T_{\mu}^{n_k}f:n\in I_k \} \rVert_{v(s)}^s \Big)^{1/s}\]
is bounded in $L^1$ for $1<s(1-\alpha-\beta)$, and \[\Big(\sum_k n_k^{\beta s} \max_{n\in I_k} |T_{\mu}^n f-T_{\mu}^{n_k}f|^s \Big)^{1/s}\]
is bounded in $L^1$ for $1-\alpha<s(1-\alpha-\beta)$.

~\hfill$\square$

\noindent {\textit Proof of Theorem \ref{thm:RittResult}.}
By Theorem 1.3. in \cite{Dunn}, there exists  a power bounded 
operator $R$ and $\gamma\in (0,1)$ 
such that 
$T=I-(I-R)^{\gamma} = R_{\nu_{\gamma}}$, where $\nu_{\gamma}$ is the probability measure on $\mathbb Z$ defined in (\ref{eqn:valpha}).  
Therefore Theorem \ref{thm:main}, Corollary \ref{cor:sup}, Propositions \ref{prop:genvar} and \ref{prop:longvar} apply to $T$.
\hfill$\square$



\begin{thebibliography}{99}

\bibitem{BC} A. Bellow and A. P. Calder\'on, {\textit A weak-type inequality for convolution products}, Harmonic analysis
and partial differential equations (Chicago, IL, 1996), 41--48, Chicago Lectures in Math., Univ.
Chicago Press, Chicago, IL, 1999.

\bibitem{BJR-power} A. Bellow, R. Jones and J. Rosenblatt, {\textit Almost everywhere convergence of powers}, G. Edgar, L. Sucheston (Eds.), Almost Everywhere Convergence, Acad. Press (1989),  99--120.

\bibitem{BJR-conv} A. Bellow, R. Jones and J. Rosenblatt, {\textit Almost everywhere convergence of convolution powers}, Ergodic
Theory Dynam. Syst. {\bf 14} (1994), 415-432.

\bibitem{Bl} S. Blunck, {\textit Analyticity and Discrete Maximal Regularity on $L^p$-Spaces}, J. Func. Anal. {\bf 183} (2001), 211--230.

\bibitem{CCL} 
G. Cohen, C. Cuny, M. Lin, {\textit Almost everywhere convergence of powers of some positive $L^p$ contractions}, J. Math. Anal. Appl. {\bf 420} (2014), 1129--1153.


\bibitem{Cuny-weak} C. Cuny, {\textit On the Ritt property and weak type maximal inequalities for convolution powers on $l^1(\mathbb Z)$}. Studia Mathematica {\bf 235} (1), (2016), 47--85.

\bibitem{Dunn} N. Dungey, {\textit Subordinated discrete semigroups of operators}, Trans.\ AMS 363 (2011), 1721--1741.

\bibitem{HH} B.H. Haak, and  M. Haase, {\textit Square Functions Estimates and Functional Calculi}, preprint (2013). \url{https://www.math.u-bordeaux.fr/~bhaak/recherche/HaHa-sqf.pdf}, \url{https://arxiv.org/abs/1311.0453}.


\bibitem{JOR} R.L. Jones, I. Ostrovskii and J. Rosenblatt, {\textit Square functions in Ergodic Theory}, Ergod.\ Th.\ \& Dynam.\ Sys. (1996), {\bf 16}, 267--305.

\bibitem{JR} R.L. Jones and K. Reinhold, {\textit Oscillation and variation inequalities for convolution powers}, Ergod.\ Th.\  \& Dynam.\ Sys.\ (2001), {\bf 21}, 1809--1829.

\bibitem{LeM-H} C. LeMerdy, $H^{\infty}$ {\textit Functional calculus and square function estimates for Ritt operators}. Rev.\ Mat.\ Iberoam.\  {\bf 30} (2014), no. 4, pp. 1149--1190. 
DOI 10.4171/RMI/811

\bibitem{LeMX-max} C. Le Merdy, and Q. Xu, {\textit Maximal theorems and square functions for analytic operators on $L^p$-spaces}, J. Lond.\, Math.\, Soc.\, (2) {\bf 86} (2012), no. 2, 343--365.

\bibitem{LeMX-Vq}C. Le Merdy, and Q. Xu, {\textit Strong q-variation inequalities for analytic semigroups}, Ann. Inst. Fourier
(Grenoble) {\bf 62} (2012), no. 6, 2069--2097 (2013).

\bibitem{LinDer} Y. Derriennic, M. Lin, {\textit Fractional Poisson equations and ergodic theorems for fractional coboundaries}, Israel J. Math. 123 (2001) 93--130.

\bibitem{Losert1} V. Losert, {\textit A remark on almost everywhere convergence of convolution powers}, Illinois J. Math. {\bf 43}
(1999), 465-479.

\bibitem{Losert2} V. Losert, {\textit The strong sweeping out property for convolution powers}, Ergodic Theory Dynam. Systems
{\bf 21} (1) (2001), 115--119.

\bibitem{Lyub}  Y. Lyubich, {\textit The single-point spectrum operators satisfying Ritt's resolvent condition}, Studia Math.\ {\bf 145} (2001), 135--142.

\bibitem{NZ} B. Nagy, J. Zemanek, {\textit A resolvent condition implying power boundness}, Studia Math {\bf 134} 2 (1999), 143--151.

\bibitem{Nev} O. Nevanlinna, {\textit Convergence of Iterations for Linear Equations}, Birkhauser, Basel, 1993.

\bibitem{Ritt} R.K. Ritt, {\textit A condition that $\lim_{n\to\infty} n^{-1} T^n = 0$}, Proc. of the AMS {\bf 4} (1953), 898-899.

\bibitem{Stein} E. Stein, {\textit On the Maximal Ergodic Theorem}, Proc. N.A.S. {\bf 47} (1961), 1894-1897.

\bibitem{Chris} C. Wedrychowicz, {\textit Almost everywhere convergence of convolution powers without finite second moment}, Ann.\, Inst.\, Fourier (Grenoble) {\bf 61} (2) (2011), 401--415.

\end{thebibliography}
\end{document}